\documentstyle{article}
\def\qed{\ifmmode\square\else\nolinebreak\hfill$\diamondsuit$\fi\par\vskip12pt}

\def\CB{{\cal B}}

\def\CE{{\cal E}}

\def\CH{{\cal H}}

\def\CM{{\cal M}}

\def\CR{{\cal R}}

\def\CV{{\cal V}}

\def\a{\alpha}
\def\b{\beta}
\def\d{\delta}
\def\D{\Delta}
\def\e{\epsilon}
\def\g{\gamma}

\def\l{\lambda}
\def\L{\Lambda}

\def\ve{\varepsilon}

\def\w{\omega}

\def\z{\zeta}

\def\lo{\left}
\def\ro{\right}
\def\pa{\paragraph}
\def\ha{{\scriptstyle{1\over2}}}
\def\be{\begin{equation}}
\def\ee{\end{equation}}

\def\st{\stackrel}
\title{\bf  Extremes of Geometric Variables with Applications to Branching
Processes.}
\author{Kosto V. Mitov\thanks{Air Force Academy G. Benkovski, Pleven,
Bulgaria. {\em email:kmitov@af-acad.bg}}\\
Anthony G. Pakes\thanks{Dept. Maths. \& Stats, Univ. of W.
Australia, 35 Stirling Highway, Crawley, W.A.
6009.{\em email:pakes@maths.uwa.edu.au}}\\
George P. Yanev\thanks{Univ. of South Florida, 140 7th Ave S.,
DAV258, St Petersburg, FL 33701, USA. {\em
email:yanev@stpt.usf.edu}}} \date{\empty}
\begin{document}
\maketitle
\begin{abstract} We obtain limit theorems for the row extrema of a
triangular array of zero-modified geometric random variables. Some
of this is used to obtain limit theorems for the maximum family
size within a generation of a simple branching
process with varying geometric offspring laws.\\
\\
{\sl Keywords \& Phrases:} Sample extrema; Geometric arrays;
Branching
processes; Varying environments; Maximum family sizes.
\end{abstract}

\section{Introduction }
\setcounter{equation}{0}
\def\theequation{\thesection.\arabic{equation}}
 It is well known (Anderson (1970)) that the geometric law is not
attracted to any max-stable law and hence maxima of independent
geometric variables cannot be approximated by a max-stable law.
Considering triangular arrays of zero-modified geometric laws
allows adjustment of the zero-class probability independently of
the success probability parameter, thus opening the possibility of
approximating row maxima and minima by simple explicit laws.
Similar results for Poisson laws are given by Kolchin et al.
(1978, \S2.6). Our motivation is closer to that of Anderson et al.
(1997) who study the Poisson and other laws. They exploit the
normal approximation to the Poisson law with large mean to show
that row maxima are approximated by the Gumbel law under certain
conditions. We obtain corresponding results which emanate from the
exponential approximation to the geometric law when its mean is
large.

For each $n=1,2,\dots$ we let $\nu_n$ be a positive integer and
$\{X_i(n):i=1,\dots,\nu_n\}$ be independent random variables with
the same zero-modified geometric law \be
P(X_i(n)=0)=1-a_n\qquad\&\qquad
P(X_i(n)=j)=a_np_n(1-p_n)^{j-1},\quad(j=1,2,\dots),\label{gl} \ee
where $0< a_n\le1$ and $0<p_n<1$. The mean for row $n$ is
$a_n/p_n$ and the distribution function is \be
F_n(x)=\cases{1-a_n(1-p_n)^{[x]} & if $x\ge0$,\cr 0 & if
$x<0$,}\label{df} \ee where $[x]$ denotes the integer part of $x$.
The standard geometric law corresponds to $a_n=1-p_n$.

In the next section we  prove limit theorems as $\nu_n\to\infty$
for the row extrema and range
$$\CM_n=\max_{1\le i\le \nu_n}X_i(n),\qquad \mu_n=\min_{1\le i\le
\nu_n}X_i(n), \quad\&\quad \CR_n=\CM_n-\mu_n,$$
 and we give examples showing our hypotheses can be satisfied. More
specifically, we find conditions which ensure the row extrema
converge in probability to infinity, and show in Theorems 1 and 3
under a further condition that normalized versions have
non-defective limit laws. Theorem 5 demonstrates their  joint weak
convergence, and Theorem 6 shows that the range is asymptotically
proportional to the maximum. Lemmas 1-4 provide context by
exhibiting possible behaviours of the extrema under differing
assumptions.

Our results for row maxima are used in \S4 to obtain corresponding
limit theorems for the maximum family size, again denoted $\CM_n$,
in the $n$-th generation of the simple branching process where the
offspring law for individuals in generation $n-1$ is the geometric
law (\ref{gl}). Of course this is precisely the case of a varying
fractional linear offspring law which has previously been studied
by Agresti (1975), Keiding and Nielsen (1975), and Fujimagari
(1980). Maxima of random variables defined on the classical
Galton-Watson tree have been studied by Arnold and Villas\~enor
(1996), Pakes (1998), and Rahimov and Yanev (1997,9). We show in
Theorems 9-11 that results from \S2 transfer to the branching
process setting through conditional limit theorems (Theorem 7
given in \S3) for the generation sizes $Z_n$. These latter results
seem to be new, and they are the strongest possible assertions,
which can be made within our restricted class of offspring laws.

\section{Behaviour of row extrema}
\setcounter{equation}{0}
\def\theequation{\thesection.\arabic{equation}}

To set our main result for $\CM_n$ in context we begin with some
elementary descriptions of its behaviour. Observe that the
distribution function of $\CM_n$ is $H_n(x):=F_n^{\nu_n}(x)$.
\pa{Lemma 1} Let $\lim_{n\to\infty}\nu_n=\infty$. (i)
$\CM_n\st{p}\to0$ if
$\nu_n a_n\to0$.\\
(ii) If $\sum \nu_na_n<\infty$ then $P(\CM_n>0 \ {\rm i.o.})=0$,
and if the rows are independent and  $\sum \nu_na_n=\infty$ then
$P(\CM_n>0 \ {\rm i.o.})=1$.\\
(iii) Let $0<\ve<1$. If the rows are independent and $
\nu_na_n>\log n+(1+\ve)\log(\log n)$ for all large $n$ then
$P(\CM_n=0 \ {\rm i.o.})=1$.  If $\nu_na_n<\log n+\log(\log n)$
then $P(\CM_n=0 \ {\rm i.o.})=0$. \pa{Proof} For (i) observe that
$H_n(x)\to1$ for all $x>0$. The remaining assertions follow from
the Borel-Cantelli lemma and elementary estimates of
$P(\CM_n>0)\leq 1-(1-a_n)^{\nu_n}$) and
$P(\CM_n=0)$.\hfill$\bigcirc$

The quantity $\a_n=\log(\nu_na_n)$ is important to our further
considerations. The following limit theorem is easily proved.
\pa{Lemma 2} Let $\lim_{n\to \infty}\nu_n=\infty$. (i) If \be
\lim_{n\to\infty}\a_n=\a\qquad(-\infty<\a<\infty)\qquad\label{2c}
\ee and \be \lim_{n\to\infty}p_n= p\label{2p} \ee then \be
\lim_{n\to\infty} H_n(x)=\cases{\exp\lo(-e^{\a}(1-p)^{[x]}\ro) &
if $x\ge0$,\cr 0 & if $x<0$.}\label{2q} \ee

(ii) If \be \lim_{n\to\infty}\a_n=\infty\qquad\&\qquad
\limsup_{n\to\infty}p_n<1,\label{2d} \ee then
$\CM_n\st{p}\to\infty$.

The limiting distribution function in (\ref{2q}) is non-defective
if $p>0$, and it is defective if $p=0$. In the latter we have
$\CM_n\st{d}\to \CM_\infty$ where
$$P(\CM_\infty=0)=1-P(\CM_\infty=\infty)=G(-\a),$$
and $G(x)=\exp(-e^{-x})$ ($-\infty<x<\infty$) is the distribution
function of the standard Gumbel law.

Theorem 1 characterizes the rate of divergence to infinity under
some further conditions. \pa{Theorem 1} Let $\lim_{n\to
\infty}\nu_n=\infty$. Assume that for some real $c$ \be
\lim_{n\to\infty} p_n=0\qquad\&\qquad \lim_{n\to\infty} \a_n
p_n=2c.\label{2a} \ee

(i) If $\lim_{n\to \infty}\a_n=\infty$,  then $c\ge 0$ and
$$p_n\CM_n-\a_n\st{d}\to \L-c,$$
where $\L$ has a standard Gumbel law.

(ii) If $\lim_{n\to\infty}\a_n=\a\quad(-\infty<\a<\infty)$, \ then
$$p_n\CM_n\st{d}\to(\L+\a)^+.$$

 \pa{Proof} Define
 $$x_n=(x+\a_n)/p_n
 \quad\&\quad
\eta_n(x)=\nu_na_n(1-p_n)^{[x_n]}\qquad(-\infty<x<\infty).$$
By
using the expansion $[x_n]=(x+\a_n)/p_n-\d_n$, where $0\le\d_n<1$
is the fractional part of $x_n$, it is easily seen that
$$\log\eta_n(x)=\a_n+[x_n]\log(1-p_n)=
\a_n+\lo({x+\a_n\over p_n}\ro)(-p_n-\ha
p_n^2+O(p_n^3))=-x+O(p_n)-\ha \a_np_n(1+o(1)).$$ The right-hand
side converges iff (\ref{2a}) holds, and the limit is $-x-c$. Thus
(\ref{2a}) is equivalent to \be
\lim_{n\to\infty}\eta_n(x)=e^{-x-c}\qquad(-\infty<x<\infty).\label{2b}
\ee
The distribution function $G_n(x)$ of $p_n\CM_n-\a_n$ equals
$[1-\nu_n^{-1}\eta_n(x)]^{\nu_n}$ if $x_n\ge0$, and equals zero
otherwise. It follows that \be
\lim_{n\to\infty}G_n(x)=\cases{G(x+c) & if $\liminf_{n\to\infty}
x_n\ge0$,\cr 0 & if $\limsup_{n\to\infty} x_n<0$.}\label{2e} \ee
 But $x_n\to\infty$ for all real $x$ iff $\a_n\to\infty$, and assertion
(i) follows. If (\ref{2c}) and (\ref{2a}) hold then $c=0$ and
$x_n\to\pm\infty$ according as $x>-\a$ or $x<-\a$, respectively.
It follows from (\ref{2e}) that
$$\lim_{n\to\infty}G_n(x)=\cases{G(x) & if $ x> -\a$,\cr
0 & if $x<-\a$,}$$ and hence $p_n\CM_n-\a\st{d}\to\max(-\a,\L)$,
and (ii) follows. \hfill$\bigcirc$

Observe that the first member of (\ref{2a}) implies that if
$n\gg1$ then $F_n(x/p_n)\approx 1-a_n+a_n(1-e^{-x})$, the
exponential approximation mentioned in \S1. Also note that Lemma
2(ii) holds under (\ref{2a}).

The assumptions of Theorem 1 can be realized. Let
$a_n=e^{-\a}\nu_n^{-\d}$ ($\d\ge0)$ and $p_n=\g \nu_n^{-\z}$
$(\z>0)$. Then (\ref{2c}) holds iff $0\le\d<1$, and then
(\ref{2a}) holds with $c=0$. The condition (\ref{2d}) holds if
$\d=1$, and then we can admit any $p_n\to0$. Now let $\d<1$ and
choose $p_n\sim A(\log \nu_n)^{-\z}$ where $A$ is a positive
constant. Then (\ref{2c}) still holds and
$$\a_np_n=A\frac{(1-\delta)\log\nu_n-c}{(\log \nu_n)^\z}\to
\cases{0 & if $\z>1$,\cr (1-\d)A & if $\z=1$,\cr \infty & if
$\z<1$.} $$ Then Theorem 1(i) holds with $c=0$ if $\z>1$ and with
$c= (1-\d)A/2$ if $\z=1$. The case $\z<1$ is an instance of
(\ref{2a}) where $c=\infty$. In this case Theorem 1(i) suggests
that the limit law is concentrated at $-\infty$, and indeed this
is true. In fact there is no affine transformation of $\CM_n$
which has a non-degenerate limit law. However we have the
following large deviation estimate,
$$\lim_{n\to\infty}{\log P(\CM_n>x_n)\over \a_np_n}=-\frac{1}{2}.$$

Further consideration of $\eta_n(x)$, defined in the proof of
Theorem 1, shows that if (\ref{2c}) holds and if (\ref{2p}) holds
with $0<p<1$ then no affine transformation of $\CM_n$ has a
non-defective limit law. The following result, generalizing the
direct assertion of Theorem 2 in Anderson (1970) and with more
explicit centering constants, shows that it is possible to
stabilize the law of $\CM_n$. The proof is similar to that for
(\ref{2b}). \pa{Theorem 2.} Set $C_n=-\a_n/\log(1-p_n)$. If
$\lim_{n\to \infty}\a_n=\infty$ and (\ref{2p}) holds with $0<p<1$,
then
$$G(\g (x-1))=\liminf_{n\to\infty}P(\CM_n-C_n\le x)\leq
\limsup_{n\to\infty}P(\CM_n-C_n\le x)=G(\g
x),\quad(-\infty<x<\infty),$$ where $\g=-\log(1-p)$.

Parallel to Lemmas 1 and 2 we have the following results for the
row minimum $\mu_n$, and they are easy consequences of its
distribution function
$$K_n(y)=\cases{1-\lo(a_n(1-p_n)^{[y]}\ro)^{\nu_n} & if $y\ge0$,\cr
0 & if $y<0$.} $$ \pa{Lemma 3} Suppose that $\nu_n\to\infty$. (i)
$\mu_n\st{p}\to0$ iff
$a_n^{\nu_n}\to0$.\\
(ii) If $\sum a_n^{\nu_n}<\infty$ then $P(\mu_n>0 \ {\rm
i.o.})=0$, and if the rows are independent and $\sum
a_n^{\nu_n}=\infty$
then $P(\mu_n>0 \ {\rm i.o.})=1$.\\
(iii) If the rows are independent and $\sum
(1-a_n^{\nu_n})=\infty$ then $P(\mu_n=0 \ {\rm i.o.})=1$, and if
$\sum (1-a_n^{\nu_n})<\infty$ then $P(\mu_n=0 \ {\rm i.o.})=0$.

\pa{Lemma 4} Let $\lim_{n\to \infty}\nu_n=\infty$. Assume that \be
\lim_{n\to\infty}\nu_n(1-a_n)=\b\qquad\&\qquad
\lim_{n\to\infty}\nu_np_n=\rho\label{2az}. \ee

(i) If $\b+\rho<\infty$, then
$$\lim_{n\to\infty}K_n(y)=\cases{1-e^{-\b-\rho[y]} & if $y\ge0$,\cr
0 & if $y<0$.} $$

(ii) If $\b=\rho=0$, then $\mu_n\st{p}\to\infty$.

The limit law in Lemma 4(i) is non-defective if $\rho>0$ and we
see that it defines a zero-modified geometric law. If $\rho=0$
then $\mu_n\st{d}\to \mu_\infty$ where
$$P(\mu_\infty=\infty)=1-P(\mu_\infty=0)=e^{-\b}.$$
Our principle result shows that if (\ref{2az}) holds with $\rho=0$
then $\mu_n$ can be centered and scaled to give a non-degenerate
limit law.   Set
$$\b_n=-\nu_n\log a_n. 
$$ \pa{Theorem 3} Let $\lim_{n\to \infty}\nu_n=\infty$. Assume that for some $b\in[0,\infty)$
\be \lim_{n\to\infty}\nu_np_n=0\qquad\&\qquad
\lim_{n\to\infty}\b_np_n=2b\label{2mb} \ee
and
 \be \lim_{n\to\infty}\b_n=\b\in[0,\infty),\label{2md} \ee then
\be \nu_np_n\mu_n\st{d}\to (\CE-\b)^+\label{2me} \ee where $\CE$
has a standard exponential law.
\pa{Proof} Let
$$y_n=(y-\b_n)/\nu_np_n \quad\&\quad \psi_n(y)=
\lo(a_n(1-p_n)^{[y_n]}\ro)^{\nu_n} \qquad (-\infty<y<\infty).$$
Observing that
$$\log\psi_n(y)=-\b_n+\nu_n\lo({y-\b_n\over \nu_np_n}-\d_n\ro)(-p_n-\ha
p_n^2+O(p_n^3)) =-y+\ha \b_np_n+\nu_n\d_np_n(1+o(1)),$$ where
$\d_n$ is the fractional part of $y_n$, it is clear that
(\ref{2mb}) is equivalent to \be
\lim_{n\to\infty}\psi_n(y)=e^{b-y},\qquad(-\infty<y<\infty).\label{2mc}
\ee
 If (\ref{2md}) holds then $b=0$ and
$y_n\to\pm\infty$ according as $y>\b$ or $y<\b$, respectively.
Since $K_n(y_n)=1-\psi_n(y)$ if $y_n>0$, $=0$ otherwise, it
follows that $\nu_np_n\mu_n+\b_n\st{d}\to\max(\b,\CE)$, and
(\ref{2me}) follows. \hfill$\bigcirc$\\

The proof shows that $\b<\infty$  is a necessary condition for a
non-defective limit law. The limit assertion (\ref{2me}) extends
for $\lim_{n\to \infty} \nu_np_n>0$ in a manner similar to Theorem
2 as follows. \pa{Theorem 4} If $\lim_{n\to\infty}
\nu_np_n=\xi\in(0,\infty)$ and (\ref{2md}) holds then
$$1-e^{-y-\b-\xi}\le \liminf_{n\to\infty} P(\nu_np_n\mu_n\le y)\le
\limsup_{n\to\infty} P(\nu_np_n\mu_n\le y) \le 1-e^{-y-\b}.$$

Observe again that (\ref{2md}) implies the second member of
(\ref{2mb}) with $b=0$, and that $a_n\to1$. Hence our assumptions
for (\ref{2md}) are precisely (\ref{2az}) with $\rho=0$. In
addition, $p_n\log \nu_n\to0$ if $\rho=0$ and $\a_n=\log
\nu_n+o(1)$. Thus, the assumptions for Theorem 1(i) are satisfied
with $c=0$ and hence we have
$$p_n\CM_n-\log \nu_n\st{d}\to \L\qquad\&\qquad
\nu_np_n\mu_n\st{d}\to(\CE-\b)^+.$$ Our next result extends this
pair of weak limit statements to joint convergence, showing in
particular that $\mu_n$ and $\CM_n$ are asymptotically
independent.
\pa{Theorem 5} If $\lim_{n\to \infty}\nu_np_n=0$ and
(\ref{2md}) hold, then
$$(p_n\CM_n-\log \nu_n,\nu_np_n\mu_n)\st{d}\to(\L,(\CE-\b)^+).$$
\pa{Proof} If $x_n=(x+\log \nu_n)/p_n$ and $y_n=y/\nu_np_n$, where
$x,y$ are real, then $x_n-y_n=(\nu_n\log
\nu_n+\nu_nx-y)/\nu_np_n\to\infty$.   Consequently for any real
$x,y$ and $n$  large enough we have
\begin{eqnarray*}\D_n(x,y):&=&P(\CM_n\le x_n,\mu_n\le
y_n)=H_n(x_n)-P(y_n<\mu_n,\CM_n\le x_n)\\
&=& H_n(x_n)\lo[1-\lo(1-{F_n(y_n)\over F_n(x_n)}\ro)^{\nu_n}\ro].
\end{eqnarray*}
The proofs of Theorems 1 and 3 show that $H_n(x_n)\to G(x)$, and
hence $F_n(x_n)\to1$, and $\nu_nF_n(y_n)\to y+\b$ if $y>0$, $\to0$
otherwise. It follows that $\D_n(x,y)\to G(x)(1-e^{-y-\b})$ if
$y>0$, $\to0$ if $y<0$.\hfill$\bigcirc$\\

Our last result shows that the row ranges $\CR_n=\CM_n-\mu_n$ are
determined
by the row maxima.\\
{\bf Theorem 6} If the conditions of Theorem 1 hold, then $\CR_n$
has the same limit behaviour as $\CM_n$. \pa{Proof} The
assumptions imply that $p_n\to0$. For any $y>0$ we have
$$-\log P(p_n\mu_n>y)=-\b_n-\nu_ny+O(\nu_np_n)\to-\infty,$$
i.e., $p_n\mu_n\st{p}\to0$. The assertion follows from Slutsky's
lemma.\hfill$\bigcirc$

\section{The simple branching process}
\setcounter{equation}{0}
\def\theequation{\thesection.\arabic{equation}}

Let $(Z_n:n\ge 0)$ denote the generation sizes of the simple
branching process with varying geometric environments,
$$Z_n=\sum_{i=1}^{Z_{n-1}} X_i(n)\qquad(n=1,2,\dots)$$
where $Z_0=1$ and the $X_i(n)$ $(i,n\ge1)$ have the same geometric
laws as in \S1, and they are mutually independent. Thus $X_1(n)$
is a generic family size of a parent in generation $n-1$, and its
probability generating function (pgf) is \be
f_n(s)={(1-s)R_n+s\over (1-s)r_n+1}\label{3a} \ee where
$$r_n=p_n^{-1}-1,\qquad R_n=p_n^{-1}-m_n,$$
and
$$m_n=f'_n(1)=a_n/p_n$$
is the mean $n$-th generation family size.

As is well known (Harris (1963)) the pgf of $Z_n$ is obtained by
functional composition,
$$\phi_n(s):=E\lo(s^{Z_n}|Z_0=1\ro)=\phi_{n-1}(f_n(s)).$$
The group structure of Mobius transformations (\ref{3a}) permits
the explicit determination
$$\phi_n(s)={(1-s)A_n+s\over (1-s)B_n+1}$$
where
$$A_n=M_n\sum_{j=1}^n R_j/M_j,\qquad B_n=M_n\sum_{j=1}^n r_j/M_j,$$
$M_0=1$ and
$$M_n=\prod_{j=1}^n m_j=E(Z_n|Z_0=1),\qquad(n\ge1).$$
The proof is (barely) indicated by Agresti (1975), and with
differing
notation.\\

The following result is fundamental.\\
{\bf Theorem 7} The conditioned process $(Z_n|Z_n>0)$ has a limit
law iff $\lim_{n\to \infty}B_n=B\in [0,\infty]$. Suppose this
condition holds. (i) If $0\le B<\infty$, then
$$(Z_n|Z_n>0)\st{d}\to \Xi$$
where
$$E(s^\Xi)={s\over 1+B-Bs}.$$
(ii) If $B=\infty$ then
$$(Z_n/B_n|Z_n>0)\st{d}\to \CE,$$
where $\CE$ has a standard exponential law.
\pa{Proof} Observe
that
$$1-\phi_n(s)={(1+B_n-A_n)(1-s)\over (1-s)B_n+1},$$
whence
$$E\lo(\lo. s^{Z_n}|\ro. Z_n>0\ro)=1-{1-\phi_n(s)\over 1-\phi_n(0)}={s\over
(1-s)B_n+1}.$$ Assertions (i) and (ii) are obvious if $B$ exists,
and the only other possibility is that $(B_n)$ has many limit
points, and hence there exists no limit law, with or without
normalization.\hfill$\bigcirc$\\

Let $T=\inf\{n:Z_n=0\}$ denote the time to extinction.\\
{\bf Theorem 8} (i) $Q:=\lim_{n\to\infty}P(Z_n=0)=1$ iff
$\lim_{n\to \infty}M_n=0$ and/or
$\sum_{j\ge1}r_j/M_j=\infty$.\\
(ii) If $Q$ exists and $B=\infty$, then \be Z_n/M_n\st{a.s.}\to
I+(1-I)\CE\label{3as} \ee
where $P(I=1)=Q=1-P(I=0)$, and $I$ and $\CE$ are independent.\\
(iii) If $Q=1$, then $P(T<\infty)=1$.
\pa{Proof} Since
$r_n-R_n=m_n-1$, we obtain
$$\sum_{j=1}^n{r_j-R_j\over M_j}=
\sum_{j=1}^n\lo(M_{j-1}^{-1}-M_j^{-1}\ro)=1-M_n^{-1}\quad\&\quad
1+B_n-A_n=M_n.$$ It follows that \be P(Z_n>0)={M_n\over
1+B_n}=\lo[M_n^{-1}+\sum_{j=1}^n r_j/M_j\ro]^{-1}.\label{3b} \ee
Assertion (i) follows. Assertion (ii) follows from Theorem 7(ii)
and the fact that $(Z_n/M_n)$ is a positive martingale which thus
converges a.s. with no further assumptions. Assertion (iii)
follows since $Z_n\st{p}\to 0$ if $Q=1$, and any a.s. convergent
subsequence eventually hits zero.
\pa{Definition 1} We say that
the environments are weakly varying if $M=\lim_{n\to\infty} M_n$
exists. If this is the case then we
have a classification similar to the trichotomy. The environments are:\\
(a) Supercritical if $M=\infty$, i.e. if $\sum_n(m_n-1)=\infty$; \\
(b) Critical if $0<M<\infty$, i.e. if $\sum_n(m_n-1)$ converges; and\\
(c) Subcritical if $M=0$, i.e. if $\sum_n(m_n-1)=-\infty$.\\

Keiding and Nielsen (1975, Theorem 2.2) prove (\ref{3as}) in the
supercritical case. This follows from Theorem 8(ii) since
$B=\infty$ if $M=\infty$. For example we can choose, as in the
classical case, $M_n\sim cm^n$,  where $c>0$ and $m>1$ are
constants and $\sup_{j\le n}r_j=O[(m/(1+\e))^n]$. Clearly $Q<1$,
Theorem 8(ii) holds, and these conditions can be realized.
However, it is possible that $Q=1$ in the supercritical case. This
occurs in our example if $p_n\to0$ sufficiently fast.
Specifically, if
$$p_n=a/b_nm^n\qquad\&\qquad a_n=mp_n$$
where $a>0$ and $\{b_n\}$ is a positive sequence, then
$r_n/M_n\sim b_n/ac$. Now choose $\sum_nb_n=\infty$. Indeed, if
$S_n$ denotes the sum in (\ref{3b}) then $P(Z_n>0)\sim S_n^{-1}$,
and we see that the probability of non-extinction can tend to zero
arbitrarily slowly. In this case (\ref{3as}) holds with $I=1$, but
Theorem 7(ii) shows that conditioning on non-extinction gives a
non-degenerate limit law.

The non-classical growth rate $M_n\sim n^\d$ is achieved by
choosing $a_n$ and $p_n$ so that $m_n=(1+n^{-1})^\d$ $(\d>0)$, and
then  $\{S_n\}$ can have a finite or infinite limit. If $\d\le1$
and $\{p_n\}$ is non-decreasing then $Q<1$ only if $p_n\to1$.

It is possible to have $B=\infty$ when $M<\infty$, but in such a
case (\ref{3b}) shows that $Q=1$, and then $Z_n/M_n\st{a.s.}\to0$.
In the critical case this occurs  iff $\sum_n r_n=\infty$, that is
$Q=1$ iff $B=\infty$ iff $\sum_n r_n=\infty$. In particular
$\sum_n r_n<\infty$ yields a non-classical restrained growth
r\'egime in the sense that $Z_n\st{d}\to I+(1-I)\Xi$.

Finally, it follows from (\ref{3b}) that $Q=1$ in the subcritical
case, and Theorem 8(iii) applies.  However it still is possible
that $B=\infty$, with Theorem 7(ii) holding. To see this, let
$m<1$ in the classical scenario above, and observe that
$$B_n\sim m^n\sum_{j=1}^n m_j r_j=\sum_{i=0}m^ir_{n-i}.$$
Hence, if $r_n\to r<\infty$ then $B=r/(1-m)$. On the other hand,
$r_n\to\infty$ if $p_n\to0$ and then Fatou's lemma shows that
$B=\infty$.

Choosing $m_n=(1+n^{-1})^{-\d}$ gives the non-classical decay
$M_n\sim n^{-\d}$, and then $B_n\sim n^{-\d}\sum_{j=1}^nj^\d r_j$.
If $r_n=n^\a l(n)$ where $\a$ is a real constant and $l(x)$ is
slowly varying, then
$$B_n\sim\cases{0 & if $\a<-1$,\cr
l(n) & if $\a=-1$, \cr \infty& if $\a>-1$.} $$
Hence $B=\infty$ if $-1<\a<0$, in which case $p_n\to1$. \\

\section{Maximum family size}
\setcounter{equation}{0}
\def\theequation{\thesection.\arabic{equation}}

Consider now the maximum $n$-th generation family size
$\CM_n:=\max_{1\le i\le Z_{n-1}}X_i(n)$. Since $\CM_n=0$ if
$Z_{n-1}=0$ we consider the conditional distribution function \be
\CH_n(x):=P(\CM_n\le
x|Z_{n-1}>0)=E\lo[F_n^{Z_{n-1}}(x)|Z_{n-1}\ro],\label{fr} \ee
where $F_n$ is defined at (\ref{df}).  Analogues of both Lemmas 1
and 2 can be given, but here we consider only the one of Lemma~2.
\pa{Theorem 9} Suppose that $\lim_{n\to \infty}p_n=p$ with
$0<p<1$. \\
(i) If $\lim_{n\to \infty}a_n= a$ and $\lim_{n\to
\infty}B_n=B<\infty$, then
$$\lim_{n\to\infty}\CH_n(x)=\frac{1-a(1-p)^{[x]}}{1+aB(1-p)^{[x]}},\qquad(x\ge 0).$$
(ii) If $\lim_{n\to \infty}B_n=\infty$ and  $\lim_{n\to
\infty}\log(a_nB_n)=\alpha $,
 then
$$\lim_{n\to\infty}\CH_n(x)=\lo[1+e^\a(1-p)^{[x]}\ro]^{-1},\qquad(x\ge 0).$$
\pa{Proof} For (i) observe that $F_n(x)\to 1-a(1-p)^{[x]}$. The
assertion follows from Theorem 7(i), (\ref{fr}), and the uniform
convergence property of the continuity theorem for probability
generating functions. For (ii), observe that the assertion of
Theorem 7(ii) is equivalent to the limit statement
$$\lim_{n\to\infty} E\lo.\lo(s^{Z_{n-1}/B_{n-1}}\ro|Z_{n-1}>0\ro)
=\lo[1+\log s^{-1}\ro]^{-1},$$ and the convergence is uniform with
respect to $s$ in the interval $[0<s'\le s\le1]$. The assertion
follows by setting $s=s_n:=F_n^{[B_{n-1}]}(x)$ and seeing that
$\{s_n\}$ has the limit (\ref{2q}).\hfill$\bigcirc$

Our next result extends Theorem 1 to the branching process
setting. The proof shows that the role played by $\nu_n$ in \S2 is
here played by $[B_{n-1}]$ and indeed that its fractional part can
be ignored. Accordingly we define $\a_n^\ast=\log(a_nB_{n-1})$,
and $\CV$ denotes a random variable having the standard logistic
distribution function $L(x)=\lo(1-e^{-x}\ro)^{-1}$, all real $x$.
\pa{Theorem 10} Suppose that $\lim_{n\to \infty}B_n=\infty$ and
for some real $c$ \be \lim_{n\to\infty} p_n=0\qquad\&\qquad
\lim_{n\to\infty} \a_n^\ast p_n=2c.\label{2aa} \ee

 (i) If
$\lim_{n\to \infty}\alpha_n^\ast=\infty$, then
$$\lo(p_n\CM_n-\a_n^\ast|Z_{n-1}>0\ro)\st{d}\to \CV-c.$$

(ii) If $\lim_{n\to \infty}\a_n^\ast=\a \qquad
(-\infty<\a<\infty)$, then
$$\lo(p_n\CM_n|Z_{n-1}>0\ro)\st{d}\to (\CV+\a)^+.$$
\pa{Proof} The proof is essentially the same as for Theorem 9(ii),
letting $s_n:=F_n^{[B_{n-1}]}(x_n)$ where $x_n=(x+\a_n^\ast)/p_n$,
just as for Theorem 1. The limit of $\{s_n\}$ is given by
(\ref{2e}) for (i), and for (ii) by $G(x)$ if $x>-\a$, $=0$
otherwise. \hfill$\bigcirc$\\

Consideration of the proofs of Theorems 1 and 9 shows that the
conditions listed in Theorem 10 are necessary for its conclusions.
The simple approach we use for Theorems 9(ii) and 10 give an
obvious analogue for Theorem 2.
\pa{Theorem 11} If the notation
and assumptions of Theorem 2 stand with $\a_n$ replaced by
$\a_n^\ast$, then
$$L(\g (x-1))=\liminf_{n\to\infty}P(\CM_n-C_n\le x)\leq
\limsup_{n\to\infty}P(\CM_n-C_n\le x)=L(\g
x),\quad(-\infty<x<\infty).$$

The various conditions in Theorems 9-11 can be satisfied, but we
will show that all but one set is satisfied by the branching
process obtained from sampling the linear birth and death process
$(\CB_t)$ at irregular times, an example mentioned due to Keiding
and Nielsen (1975). Let $Z_n=\CB_{t_n}$ where $0<t_n<t_{n+1}\to
t_\infty\le\infty$. If $\l$ and $\mu$ are the birth and death
rates, respectively, and $d_n=t_n-t_{n-1}$, then $a_n=m_np_n$,
$$p_n=\cases{{\l-\mu\over \l m_n-\mu} & if $\l\not=\mu$,\cr
\lo(1+\l d_n\ro)^{-1} & if $\l=\mu$,} \quad\&\quad
m_n=e^{(\l-\mu)d_n}.$$  and
$$B_n=\cases{{\l(M_n-1)\over \l-\mu} & if $\l\not=\mu$,\cr
\l t_n & if $\l=\mu$,} \quad\&\quad M_n=e^{(\l-\mu)t_n}.$$ The
environments are weakly varying in this case, and our criticality
classification coincides with the standard one if
$t_\infty=\infty$, and they are critical if $t_\infty<\infty$. The
latter case is vacuous as far as Theorems 9-11 are concerned,
since $B<\infty$ and $p=1$, i.e., there exists no non-degenerate
conditional limit law for $\CM_n$. So we assume that
$t_\infty=\infty$ and that $d_n\to d\le\infty$.

If the environments are subcritical, $\l<\mu$, then $B<\infty$,
$p_n\to p$ with
$$0<p={\mu-\l\over \mu-\l m_\infty}<1\qquad\&\qquad 0\le
m_\infty=e^{-(\mu-\l)d}<1$$ and $a_n\to a=pm_\infty $. Hence
Theorem 9(i) holds and the limit law is degenerate at the origin
if $d=\infty$, but not otherwise.

If $\l\ge\mu$ then $B=\infty$, but Theorem 9(ii) cannot hold
because $a=(1-\l/\mu)/(1-\mu/m_\infty)>0$, whence $\a_n^\ast\to
\a<\infty$ is violated. We show that the conditions of Theorems 10
and 11 can be satisfied. First, suppose that $d=\infty$, in which
case $p=0$. In the supercritical case $\l>\mu$ we further suppose
that the second member of (\ref{2aa}) is satisfied, i.e., \be
\lim_{n\to\infty}(t_n/m_n)={2c\over
\l-\mu}\in[0,\infty).\label{3cc} \ee For example, this condition
is satisfied with $c=0$ if $d_n>\ell(n)\to\infty$ where $\ell(x)$
is a slowly varying function. If (\ref{3cc}) holds then
$$\lo((\CM_n/m_n)-(\l-\mu)t_n\ro|
Z_{n-1}>0\st{d}\to \CV-c.$$

Now suppose that $\l=\mu$ and $t_n=n^\d\ell(n)$ where $\d\ge1$,
$\ell(x)$ is slowly varying and that $t_n/n\to\infty$. Then
(\ref{2aa}) holds with $c=0$, $a_nB_{n-1}\sim n$, and again
Theorem 10(i) holds in the form
$$\lo(\lo.{\CM_n\over \l\d n^{\d-1}\ell(n)}-\log n\ro|
Z_{n-1}>0\ro)\st{d}\to {\cal V}.$$ On the other hand, if
$t_{n-1}/t_n\to\tau\in(0,1)$ then
$anB_{n-1}\to\w:=\tau/(1-\tau)>1$ and $\a_n\to\a=\log\w$. Hence
Theorem 10(b) holds in the form
$$\lo(\lo.{\CM_n\over (1-\tau)\l t_n}\ro|Z_{n-1}>0\ro)\st{d}\to (\CV+\a)^+.$$
If $\l\ge\mu$, $t_\infty=\infty$ but $d<\infty$ then $0<p<1$ and
Theorem 11 holds.

We end by observing that analogues of results in \S2 for the
minimum and range can be taken into the branching process context,
but we leave this as an exercise for the reader.

\end{document}